\documentstyle[12pt,twoside,leqno]{article}
\def\bbfR{I\!\!R}
\newcommand{\rf}[1]{(\ref{#1})}

\def\B1b{{{\cal B}_{1}^{-\beta,\infty}}}
\def\H{{\cal H}}
\def\cbdu{\hfill{$\Box$}}
\def\X{{\cal X}}
\def\Nop{{\cal N}}
\newcommand{\proof}{\noindent {\bf Proof.} }
\newcommand{\rn}{{I\!\!R}^n}
\newcommand{\BR}{{I\!\!R}}
\newcommand{\dd}{\mbox{\,d}}
\newtheorem{theorem}{Theorem}[section]
\newtheorem{proposition}{Proposition}[section]
\newtheorem{lemma}{Lemma}[section]
\newtheorem{corollary}{Corollary}[section]

\newcounter{remark}
\newenvironment{remark}%
{\medskip \stepcounter{remark} \noindent {\bf Remark
\arabic{section}.\arabic{remark}.}}{\rm \cbdu}

\begin{document}
$ $
\vspace{1cm}

\centerline{\large\bf ON ZERO MASS SOLUTIONS}
\medskip

\centerline{\large\bf OF VISCOUS CONSERVATION LAWS}

\

\

\begin{center}
{\large Grzegorz Karch$^1$ and Maria Elena Schonbek$^2$}\\
 $\;$                       \\
{\small $^1$Instytut Matematyczny, Uniwersytet Wroc{\l}awski}\\
{\small pl. Grunwaldzki 2/4, 50--384 Wroc{\l}aw, Poland}\\
{\small \it karch@math.uni.wroc.pl}\\
{\small $^2$Department of Mathematics, University of California}\\
{\small Santa Cruz, CA 95064,  U.S.A.}\\
{\small \it schonbek@math.ucsc.edu}
\end{center}
\vspace{0.2cm}

\

\

\noindent{\small 
{\bf Abstract.}
 In the paper, we consider the large time behavior of solutions to
 the convection-diffusion
equation $u_t-\Delta u +\nabla\cdot f(u)=0$ in $\bbfR^n\times
[0,\infty)$, where $f(u) \sim u^q$ as $u\to 0$. 
Under
the assumption that 
$q \geq  1+1/(n+\beta)$ and 
the initial condition $u_0$
satisfies:
$u_0\in L^1(\bbfR^n)$, $\int_{\bbfR^n} u_0(x)\;dx=0$, and 
$\|e^{t\Delta}u_0\|_{L^1(\bbfR^n)}\leq Ct^{-\beta/2}$ for fixed
$\beta\in (0,1)$, all $t>0$, and a constant $C$, we show that 
the $L^1$-norm of the solution to the convection-diffusion equation
decays with the rate $t^{-\beta/2}$ as $t\to\infty$. 
Moreover, we prove that, for small initial  conditions,   the
exponent
 $q^*=1+1/(n+\beta)$ is critical  in the following sense.
For $q>q^*$ the large time  behavior 
in $L^p(\bbfR^n)$, $1\leq p\leq \infty$, 
of solutions is described by 
self-similar  solutions to the linear heat equation. For $q=q^*$,
we prove that the convection-diffusion  equation with
$f(u)=u|u|^{q^*-1}$ has a family of self-similar solutions which
play an important role in the large time asymptotics of general
solutions.

\footnote[0]{{\it 2000 Mathematics Subject Classification}: 35B40,
35K55.
}
\footnote[0]{{\it Key words and phrases}: the Cauchy problem,
the convection-diffusion equation, large time behavior of solutions,
self-similar solutions.
}
}

\newpage
\section{Introduction}
In this paper, we study the large time behavior of solutions 
$u=u(x,t)$ $(x\in\bbfR^n, \; t>0)$ to the
Cauchy problem for the nonlinear convection-diffusion equation
\begin{eqnarray}
&u_t-\Delta u+a\cdot \nabla (u|u|^{q-1})=0,&
\label{eq}\\
&u(x,0)=u_0(x),&\label{ini}
\end{eqnarray}
where $q>1$ and the vector $a\in\bbfR^n $ are fixed.
The assumptions $u_0\in L^1(\bbfR^n)$ and
$\int_{\bbfR^n}u_0(x)\;dx =0$ will also be required.

The typical nonlinear term occurring in hydrodynamics in the
one-di\-men\-sion\-al case has the form $u u_x= (u^2/2)_x$ (as in the
case of the viscous Burgers equation). The most obvious
generalization of this nonlinearity consists in replacing the square
by a power $u^q$ where $q$ is a positive integer. 
Here, however, we intend  to observe a more subtle
interaction of the nonlinearity with dissipation,
consequently, we need to consider a continuous range of  parameters
$q$. The problem then appears with the definition of $u^q$ for
negative $u$ and for non-integer $q$.
In order to avoid this difficulty, we chose the nonlinear
term of the from $a\cdot \nabla(u|u|^{q-1})$. This was done only
to shorten notation in this report. 
Note that, in fact, the
following property of the nonlinearity will  be  essential
throughout this work: 
\begin{itemize}
\item 
the nonlinear term in $\rf{eq}$
has the form $\nabla\cdot f(u)$ where the $C^1$-vector function $f$
satisfies
$
|f(u)|\leq C|u|^q,
$
$
|f'(u)|\leq C|u|^{q-1}
$
for every $u\in\bbfR$, $q>1$, and a constant $C$.
Moreover, if the balanced case is considered (i.e.
$q=1+1/(n+\beta)$), the limits 
$$
\lim_{u\to 0^-} f(u)/|u|^q,
\quad \mbox{and} \quad
\lim_{u\to 0^+} f(u)/|u|^q 
$$
should exist  and  the both should be different from 0.
\end{itemize}

Recent publications developed versatile functional analytic tools to
study  the long time behavior of solutions of this initial value
problem. 

Concerning the decay  of solutions of \rf{eq}-\rf{ini} and, more
generally, of
scalar parabolic conservation laws of the form
$
u_t -\Delta u+ \nabla\cdot f(u) = 0
$
with integrable initial conditions,  
Schonbek \cite{S} was the first who proved that the $L^2$-norm tends
to 0 as $t\to \infty$ with the  rate $t^{-n/4}$. To deal with this
problem, she introduced the so-called {\it Fourier splitting method.}
The results from  \cite{S} were extended in the later work \cite{S1}, 
where the decay of solutions in $L^p(\bbfR^n)$, 
($1\leq p\leq \infty$)  was obtained, again,  
by a method based on the Fourier splitting technique.
It was emphasized in \cite{S1} that
the decay rates  are the same as
for the underlying linear equations. 

Next, Escobedo and Zuazua \cite{EZ} proved  decay estimates of the
$L^p$-norms of solutions by a different method under more general
assumptions on nonlinearity and under less restrictive assumptions on
initial data. Finally, by the use of the {\it logarithmic Sobolev
inequality},  Carlen and Loss \cite{CL} showed that solutions 
of viscous conservation laws satisfy
$$
\|u(\cdot,t)\|_p\leq Ct^{-(n/2)(1/r-1/p)}\|u_0\|_r
$$
for each $1\leq r\leq p\leq \infty$, all $t>0$, and a numerical
constant $C>0$ depending on $p$ and $q$, only.
Here, we would also like to recall results on algebraic decay rates
 of solution to 
systems of parabolic
conservation laws, obtained by Kawashima \cite{KA}, Hopf and Zumbrun
\cite{HZ}, 
Jeffrey and Zhao \cite{JZ}, and Schonbek and S\"uli \cite{SS}. 
Smallness assumptions on initial conditions were often imposed in
those papers.

The first term of the asymptotic expansion 
was studied as the next step in  analysis of the long time behavior
of solutions to \rf{eq}-\rf{ini}.
Assuming that $u_0\in L^1(\bbfR^n)$, roughly speaking, these results,
cf. e.g. \cite{CTPL,KA,EZ,EVZ1,EVZ2,DuZ,DuZ2, BKW1,BKW2, BKW3,K3,K4},
fall into three cases:

\begin{itemize}
\item 
{\it Case I:}
$q>1+1/n$, when the asymptotics is linear, i.e. 
 \begin{equation}
 t^{(n/2)(1-1/p)}\|u(\cdot,t)-MG(\cdot,t)\|_p\to 0\ \mbox{as\
 }t\to\infty,
 \label{lin}
 \end{equation}
where $M=\int_{\bbfR^n} u_0(x)\, dx$, 
$G(x,t)=(4\pi t)^{-n/2}\exp(-|x|^2/(4t))$ 
is the fundamental solution of  the heat equation. 
Hence, this case can be classified as {\it weakly nonlinear}, since
in this situation the linear diffusion prevails and the nonlinearity 
is asymptotically negligible.

\item
{Case II:}
 $q=1+1/n$, when 
 \begin{equation}
 t^{(n/2)(1-1/p)}\|u(\cdot,t)-U_M(\cdot,t)\|_p\to 0\ \mbox{as\
 }t\to\infty,
 \label{auto}
 \end{equation}
where $U_M(x,t)=t^{-n/2}U_M(xt^{-1/2},1)$ is the self-similar
solution of 
\rf{eq} with $u_0(x)=M\delta_0$. Here, diffusion and 
the convection are balanced, and the asymptotics is determined by 
a special solution of a~nonlinear equation.

\item
{Case III:}
  $1<q<1+1/n$, when the convection points in the $x_n$-direction
  (i.e. $a=(0,...,0,1)$).
  Here
 \begin{equation}
 t^{(n+1)(1-1/p)/(2q)}\|u(\cdot,t)-U_M(\cdot,t)\|_p\to 0\ \mbox{as\
 }t\to\infty,
 \label{hyp}
 \end{equation}
holds, 
where $U_M$ is a particular self-similar solution of the partly
viscous 
conservation law 
$U_t-\Delta_yU+{\partial\over{\partial x_n}}(U|U|^{q-1})=0$ such that 
$u_0(x)=M\delta_0$ in the sense of measures. Here $x=(y,x_n)$, 
$y=(x_1,\dots,x_{n-1})$, and 
$\Delta_y=\sum_{j=1}^{n-1}{\partial^2\over{\partial x_j^2}}$.
Hence, the asymptotics 
of solutions is determined by solutions of an equation with strong
convection 
and partial dissipation. 
\end{itemize}

Finally, we recall that, in the weakly nonlinear case,
Zuazua \cite{Z} found, for  solutions to
\rf{eq}-\rf{ini}, the second
 order term in the
asymptotic expansion as $t\to\infty$. He observed that  asymptotic
behavior of the solution differs depending if $q$ satisfies
 $1+1/n<q<1+2/n$, $q=1+2/n$, or $q>1+2/n$. 
Analogous results for L\'evy conservation laws were obtained in
\cite{BKW1,BKW2}, and for convection-diffusion equations with
dispersive effects in \cite{K3,K4}. Related results on the stability 
in $L^1(\bbfR^n)$ of traveling waves 
(or shock waves)  in scalar viscous conservation laws can be found 
in the papers by Serre \cite{Se1} and  Freist\"uhler and Serre \cite{Se2}.
Some results on the $L\sp 1$-stability of the zero solution of degenerate 
convection-diffusion equations can be found in the article by
Feireisl and Lauren\c cot \cite{FL}.

\bigskip

Here, we assume that $M=\int_{\bbfR^n}u_0(x)\;dx
=\int_{\bbfR^n} u(x,t)\;dx=0$, thus the corresponding self-similar 
intermediate asymptotics in \rf{lin}-\rf{hyp} are equal to 0 for
every
$q>1$. Moreover, for $p=1$ the asymptotic formulae in
\rf{lin}-\rf{hyp}
say nothing else but
$\|u(\cdot,t)\|_1\to0$  as $ t\to\infty.$

The goal of this paper is to find self-similar asymptotics in
$L^p(\bbfR^n)$ 
 of solutions to \rf{eq}-\rf{ini}
with $M=0$ imposing additional  conditions on the initial data.
We assume that $u_0$ satisfies $\|e^{t\Delta}u_0\|_1\leq
Ct^{-\beta/2}$ for some
$\beta \in (0,1)$, all $t>0$, and $C$ independent of $t$. Such a
decay estimate of 
solutions to the linear heat equation is optimal for a large class of
initial conditions 
(cf. Propositions \ref{prop-L1-lin} and \ref{lem:lin:ss} below).
Under these assumptions, we  improve the known algebraic
decay rates of the solutions to \rf{eq}-\rf{ini}
in the $L^p$-norms for every  $1 \leq p \leq \infty$.  In addition,
 if the initial data are sufficiently small,
we discover the new  critical exponent 
$q^*=1+1/(n+\beta)$
such that
\begin{itemize}
\item
for $q>q^*$ the asymptotics of solutions to \rf{eq}-\rf{ini} is
linear and described 
by  self-similar solutions to the heat equation (cf. Corollaries
\ref{cor-lin-self} and \ref{cor-Lp}, below); 

\item
 $q=q^*$ corresponds to the balanced case, and the asymptotics of
 solutions corresponding to suitable small initial conditions is
 described
 by a new class of self-similar solutions to the nonlinear equation
 \rf{eq} (cf. Theorem~\ref{asymp} and the discussion in Section 5.). 
\end{itemize}

In the next section of this paper, we briefly present the main
results.
The  results corresponding to the 
 case $q>q^*$ are contained in Section 3. Section~4 considers the  
case when the exponent  $q=q^*$  is  critical.
In Section 5, we explain how to derive, from our general theorems,
self-similar solutions to the nonlinear equation \rf{eq} and how to
study large time asymptotics of general solutions.
In the last section, we discuss possible applications of our ideas to
other
equations such 
as the  Navier-Stokes equations, the KdV-Burgers equation, and the BBM-Burgers equation.

\medskip

{\bf Notation.}
The notation to be used is mostly standard. For
$1\leq p\leq \infty$, the $L^p$-norm   of a Lebesgue
measurable real-valued function defined on $\bbfR^n$
is denoted by $\|v\|_p$. We will always  denote by
$\|\cdot\|_{\cal X}$ the norm of any other Banach space $\cal X$
used in this paper.

If $k$ is a nonnegative integer,
$W^{k,p}(\bbfR^n)$ will be the Sobolev space consisting of
functions in $L^p(\bbfR^n)$ whose generalized derivatives up to order
$k$ belong  to $L^p(\bbfR^n)$.

The Fourier transform of $v$ is defined as  
$\widehat v(\xi)\equiv (2\pi)^{-n/2}\int_{\bbfR^n}
e^{-ix\xi} v(x)\;dx$.

Given a multi-index $\gamma=(\gamma_1, ..., \gamma_n)$, we denote 
$\partial^\gamma =\partial^{|\gamma|}/\partial_{x_1}^{\gamma_1}...
\partial_{x_n}^{\gamma_n}$. On the other hand, for $\beta>0$, the
operator $D^\beta$ is defined via the Fourier transform as 
$\widehat{(D^\beta w)}(\xi)=|\xi|^\beta|\widehat w (\xi)$.

The letter $C$  will denote generic positive constants, which do not
depend on $t$ and may vary from line to line during computations.

\setcounter{equation}{0}
\section{Main  results and comments}

We recall that for every $u_0\in L^1(\bbfR^n)$, the Cauchy
problem 
\rf{eq}-\rf{ini} has a unique solution in $C([0,\infty);
L^1(\bbfR^n))$ satisfying   
$$
 u\in C((0,\infty); W^{2,p}(\bbfR^n)) \cap C^1((0, \infty),
 L^p(\bbfR^n))
$$
for all $p\in (1,\infty)$.
The proof is based  on a standard iteration procedure
involving the integral representation of solutions of
\rf{eq}-\rf{ini} 
\begin{equation}
u(t)=e^{t\Delta}u_0-\int_0^ta\cdot\nabla e^{(t-\tau)\Delta}
(u|u|^{q-1})(\tau)\,d\tau \label{duhamel}
\end{equation}
(see, e.g. \cite{EZ} for details).
Here, $e^{t\Delta}u_0$ is the solution to the linear heat equation
given by the convolution of the initial datum $u_0$
with the Gauss-Weierstrass kernel
$G(x,t)=(4\pi t)^{-n/2} \exp(-|x|^2/(4t))$.
Formula \rf{duhamel}
will be  one of the main tools used in the analysis of the
long time behavior of solutions.

Let us also recall that sufficiently regular solutions of
\rf{eq}-\rf{ini}
satisfy the estimate
\begin{equation}
\|u(\cdot,t)\|_p\leq C(p,r)t^{-(n/2)(1/r-1/p)}
\|u_0\|_r \label{CL}
\end{equation}
for all $1\leq r\leq p\leq\infty$, all $t>0$, and a constant $C(p,r)$ 
depending on $p$ and $r$, only.
Inequalities \rf{CL}
are due to Carlen and Loss \cite[Theorem 1]{CL}. We also refer the
reader to \cite{BKW1, BKW2} where counterparts of \rf{CL} were proved 
for more general equations: so-called L\'evy conservation laws.

\bigskip

Section 3 contains the analysis of the large time
asymptotics of solutions to the linear heat equation.
Easy calculations show that for every $u_0\in L^1(\bbfR^n)$
such that $\int_{\bbfR^n}u_0(x)\;dx=0$ we have $\|e^{t\Delta}
u_0\|_1\to 0$ 
as $t\to\infty$. 
The following proposition asserts the existence of a large class
of initial
conditions for which 
the large time behavior of $e^{t\Delta}u_0$ is self-similar.
Here, we need the notion of the Riesz potential $I_\beta$
and the fractional derivative $D^\beta$ 
defined in the Fourier variables as
\begin{equation}
\widehat{(I_\beta w)}(\xi)={\widehat w(\xi) \over |\xi|^\beta}
\quad \mbox{and}\quad
\widehat{(D^\beta w)}(\xi)=|\xi|^\beta {\widehat w(\xi)}.
\label{ID}
\end{equation}

\begin{proposition}
\label{prop-L1-lin}
Let $\beta >0$ and $\gamma= (\gamma_1,....,\gamma_n) $ be a
multi-index with $\gamma_i \geq 0$.
Assume that $I_\beta u_0\in L^1(\bbfR^n)$.
 Denote
\begin{equation}
A=\lim_{|\xi|\to 0} {\widehat u_0(\xi)\over |\xi|^\beta}
=\int_{\bbfR^n} (I_\beta u_0)(x)\;dx. \label{lim:A}
\end{equation}
Then
\begin{equation}
\|\partial^\gamma e^{t\Delta}u_0\|_1\leq
Ct^{-\beta/2-|\gamma|/2}\|I_\beta u_0\|_1
\label{lin-L1}
\end{equation}
for all $t>0$ and $C=C(\beta,\gamma)$ independent of $t$ and $u_0$;
moreover, 
\begin{equation}
t^{\beta/2+ |\gamma|/2} \|\partial^\gamma e^{t\Delta}u_0(\cdot)-A
\partial^\gamma D^{\beta} G(\cdot,t)\|_1\to 0
 \label{as-lin-L1}
\end{equation}
as $t\to\infty$. 
\end{proposition}

Let us emphasize that  we do not
assume that $\beta<1$ in Proposition \ref{prop-L1-lin}. 
This condition only becomes necessary when the
convection term is present. Here, we also refer the reader to
Proposition \ref{lem:lin:ss} where self-similar asymptotics of the
heat semigroup is
studied under more general assumptions on initial data.

\medskip

In our first theorem on the large time behavior of solutions to the
nonlinear problem \rf{eq}-\rf{ini}, 
we  assume the decay of
$\|e^{t\Delta} 
u_0\|_1$
with a given rate and we prove that the same decay estimate holds
true 
for solutions to  \rf{eq}-\rf{ini}.

\begin{theorem}\label{tw:L1dec}
Fix $0<\beta<1$.
Assume that  $u_0\in L^1(\bbfR^n)\cap L^q(\bbfR^n)$ 
satisfies the inequality
\begin{equation}
\|e^{t\Delta}u_0\|_1\leq Ct^{-\beta/2}
\label{dec:u0}
\end{equation}
for all $t>0$ and a constant $C$ independent of $t$.
Let $u$ be the solution to \rf{eq}-\rf{ini} with $u_0$
as the initial datum.
If $q > 1+1/n$, then there exists a constant $C$ such that 
\begin{equation} 
\|u(\cdot,t)\|_1\leq C(1+t)^{-\beta/2}\label{L1:dec}
\end{equation}
for all $t>0$.
The estimate \rf{L1:dec} holds also true for  $1+1/(n+\beta)\leq q \leq
1+1/n$, with
$0<\beta<1$,
 provided $u_0\in L^1(\bbfR^n)\cap L^\infty(\bbfR^n)$
and 
$\sup_{t>0}t^{\beta/2}\|e^{t\Delta}u_0\|_1$ is sufficiently small.
\end{theorem}

\begin{remark}
The  assumption \rf{dec:u0} means that $u_0$ belongs to the
homogeneous Besov space $\B1b$ (cf. \rf{B-norm}, below) which 
will play an important role in the analysis of the balanced case
$q=1+1/(n+\beta)$.  
\end{remark}

\medskip

The approach formulated in Theorem \ref{tw:L1dec}, saying that the 
decay estimates imposed on the heat semigroup lead to the analogous 
estimates of solutions to a nonlinear problem, appears in several
recent 
papers. Here, we would like only to recall (the list is by no mean
exhaustive)
the works on the Navier-Stokes system by Schonbek \cite{SNS}  and
Wiegner \cite{W} where the 
$L^2$-decay of solutions was studied as well as by Miyakawa \cite{M}
where  decay of the $L^1$-norm and ${\cal H}^p$-norms (the Hardy
spaces)
of weak solutions was shown. Moreover, our results extend
essentially 
the recent paper by Schonbek and S\"uli \cite{SS} where general
conservation
laws 
were considered.

\medskip

If we combine the decay from \rf{L1:dec} with inequalities \rf{CL},
we obtain the improved $L^p$-decay of solutions to \rf{eq}-\rf{ini}.
Moreover, applying such estimates to \rf{duhamel} we find the
asymptotics of solutions for $q>1+1/(n+\beta)$. The following
corollary contains these results.

\begin{corollary}
\label{cor-Lp0}
Under the assumptions of Theorem \ref{tw:L1dec}, for every $p\in
[1,\infty]$ and $\beta \in (0,1)$,
there exists $C=C(u_0,p)$ independent of $t$   
such that
\begin{equation}
\|u(\cdot,t)\|_p\leq C(1+t)^{-(n/2)(1-1/p)-\beta/2} \label{u-Lp0}
\end{equation}
for all $t>0$. Moreover, for  $q>1+1/(n+\beta)$ and for every
$p\in[1,\infty]$ it follows 
\begin{equation}
t^{(n/2)(1-1/p)+\beta/2}\|u(\cdot,t)-
e^{t\Delta}u_0(\cdot)\|_p \to 0 
\quad \mbox{as}\quad  t\to\infty.
\label{ue-Lp0} 
\end{equation} 
\end{corollary}

A slightly stronger version of this corollary is formulated and
proved in the next section (cf. Corollary \ref{cor-Lp}, below).
Here, we only emphasize that combining \rf{ue-Lp0} with Proposition
\ref{prop-L1-lin} we obtain that  
the large time 
behavior of solutions to \rf{eq}-\rf{ini} with  $ q> 1 +1/n$ (or,
if the data are sufficiently small, for $q>1+1/(n+\beta)$) is
described by special self-similar solutions to the heat equation.
 This is worth  stating more precisely.

\begin{corollary} \label{cor-lin-self}
Under the assumptions of Theorem \ref{tw:L1dec} and 
Proposition 
\ref{prop-L1-lin} (or Proposition \ref{lem:lin:ss} with
$\ell(\xi)=|\xi|^\beta$, see Section 3)
the solution to \rf{eq}-\rf{ini} with $q>1+1/(n+\beta)$ satisfies 
$$
t^{(n/2)(1-1/p)+\beta/2}\|u(\cdot,t)-A D^\beta G(\cdot,t)\|_p\to 0
\quad \mbox{as} \quad
t\to\infty.
$$
\end{corollary}

\medskip

Our next results, studied in  Section 4, correspond to the balanced
case
$$
q=q^*=1+{1\over n+\beta}
$$
for some fixed $0<\beta<1$.
We will work in the homogeneous Besov space $\B1b$
defined by
$$
\B1b=\{v\in {\cal S}'(\bbfR^n)\;:\; \|v\|_\B1b<\infty\},
$$ 
where ${\cal S}'(\bbfR^n)$ is the space of tempered distributions and
the norm 
is given by 
\begin{equation}
\|v\|_\B1b\equiv\sup_{s>0}
s^{\beta/2}\|e^{s\Delta}v\|_{1}.\label{B-norm}
\end{equation}
The standard way of defining norms in Besov spaces is based on the
Paley-Littlewood dyadic decomposition. The choice of the
equivalent norm \rf{B-norm} allows us  to simplify several
calculations. Recall here that Proposition \ref{prop-L1-lin}
describes a large subset in $\B1b$ of
initial conditions  $u_0$.

\medskip
                                              
Section 4 contains the proofs of   two main theorems. The first one 
 provides  a construction of global-in-time solutions to
\rf{eq}-\rf{ini} with $q=1+1/(n+\beta)$ and  suitably small
initial data in the space $\B1b$.
The second theorem gives asymptotic  stability  of solutions 
in the balanced case. The precise statement of the theorems is the
following.

\begin{theorem} \label{exist:global}
Fix $\beta\in (0,1)$ and put $q=1+1/(n+\beta)$.
There is $\varepsilon >0$ such that for each $u_0\in\B1b$
satisfying $\|u_0\|_\B1b <\varepsilon$ there exists a solution of
\rf{eq}-\rf{ini} for all $t\geq 0$ in the space
\begin{eqnarray*}
\X&\equiv& {\cal C}([0,\infty):\B1b)\\
&& \quad \cap \;  
\{u:(0,\infty)\to L^q(\bbfR^n)\;:\;
\sup_{t>0}t^{(n/2)(1-1/q)+\beta/2}\|u(t)\|_q< \infty\}.
\end{eqnarray*}
This is the unique solution satisfying the condition 
$$
\sup_{t>0}t^{(n/2)(1-1/q)+\beta/2}\|u(t)\|_q\leq 2 \varepsilon.$$
\end{theorem}

\medskip

\begin{theorem}\label{asymp}
Let the assumptions from Theorem \ref{exist:global} hold true.
Assume that $u$ and $v$ are two solutions of \rf{eq}-\rf{ini}
constructed in Theorem \ref{exist:global}
corresponding to the initial data $u_0,v_0\in\B1b$, respectively.
Suppose that
\begin{equation}
\lim_{t\to\infty}t^{\beta/2}\|e^{t\Delta}(u_0-v_0)\|_1=0.
\label{e:u0v0}
\end{equation}
Choosing   $\varepsilon>0$ in Theorem \ref{exist:global}
sufficiently small,  we have
\begin{equation}
\lim_{t\to\infty}t^{(n/2)(1-1/p)+\beta/2}\|u(\cdot,t)-v(\cdot,t)\|_p=
0\label{e:utvt}
\end{equation}
for every $p\in [1,\infty]$.
\end{theorem}

\bigskip

In Section 5, we show  how to use Theorem \ref{exist:global} in
order to obtain  self-similar solutions
to equation \rf{eq} with 
the critical exponent   $q= 1 +1/(n+\beta)$. 
Moreover, we explain the role of self-similar solutions in the large
time asymptotics of other solutions to \rf{eq}-\rf{ini}.

\setcounter{equation}{0}
\section{Asymptotics of solutions for $q>1+1/(n+\beta)$}

As noted in  Section 2, the first problem is to find a class of
data
 that will insure the decay of solutions to the heat equation in
 $L^1(\bbfR^n)$. This 
 is obtained in Proposition \ref{prop-L1-lin}, where this class of
 data is 
shown to be constituted  by functions such that their convolutions
with Riesz potentials 
 lie in $L^1(\bbfR^n)$. Now, we establish  Proposition
 \ref{prop-L1-lin}.

\medskip

\noindent {\bf Proof of Proposition \ref{prop-L1-lin}.}
Let us note that the limit in \rf{lim:A} exists, since 
$\widehat u_0(\xi)/|\xi|^\beta$ is  continuous  
as the Fourier transform of an integrable function $I_\beta u_0$.

First, we prove that 
$\partial^\gamma D^{\beta} G(\cdot,1)\in L^1(\bbfR^n)$. 
Obviously, $\partial^\gamma D^{\beta} G(\cdot,1)$ is bounded and
continuous because 
its Fourier transform $(i\xi)^\gamma |\xi|^\beta e^{-|\xi|^2}$ is
integrable. 
Moreover, it follows from \cite[Ch. 5, Lemma 2]{St70} that for every
$\beta>0$ there exists a finite measure $\mu_\beta$ on $\bbfR^n$ given by
$$
\widehat \mu_\beta(\xi)={|\xi|^\beta\over (1+|\xi|^2)^{\beta/2}}.
$$
Hence, $\partial^\gamma D^{\beta}
G(\cdot,1)=\mu_\beta*K_{\beta,\gamma}$ where the function 
$K_{\beta,\gamma}$ is defined via the Fourier transform as
$\widehat K_{\beta,\gamma}(\xi) =(i\xi)^\gamma (1+|\xi|^2)^{\beta/2}
e^{-|\xi|^2}$. It is easy to prove that 
$K_{\beta,\gamma}\in {\cal S}(\bbfR^n)$ (the Schwartz class of
rapidly
decreasing smooth function), and this implies the integrabilty of 
$\partial^\gamma D^{\beta} G(\cdot,1)$ for every multi-index
$\gamma$.

Now,  the change of variables yields that $\partial^\gamma
D^{\beta}
G(x,t)$
has the self-similar form:
\begin{equation}
\partial^\gamma D^{\beta} G(x,t) = t^{-n/2-\beta/2-| \gamma|/2}
(\partial^\gamma D^\beta
G)(x/\sqrt t,1)
\label{DbG-self}
\end{equation}
for all $x\in\bbfR^n$ and $t>0$.

To prove \rf{lin-L1}, use the Young inequality 
for the convolution, and thus by \rf{DbG-self} it follows
\begin{eqnarray*}
\|\partial^\gamma e^{t\Delta} u_0\|_1 &=& \|\partial^\gamma D^\beta
G(t)*I_\beta u_0\|_1\\
&\leq & \|\partial^\gamma D^\beta G(\cdot,t)\|_1 \|I_\beta u_0\|_1\\
&\leq & t^{-\beta/2- |\gamma|/2}\|\partial^\gamma D^{\beta}
G(\cdot,1)\|_1\|I_\beta u_0\|_1
\end{eqnarray*}
for all $t>0$.

For the proof of \rf{as-lin-L1},  observe that the change of
variables $z=x/\sqrt t$ combined with   \rf{DbG-self} leads to the
following
expression
\begin{eqnarray}
&&\hspace{-0.5cm}t^{\beta/2+ |\gamma|/2} \|\partial^\gamma
e^{t\Delta}u_0(\cdot)-A
\partial^\gamma D^{\beta} G(\cdot,t)\|_1\nonumber\\
&&=t^{\beta/2+ |\gamma|/2}
\int_{\bbfR^n} \left|\int_{\bbfR^n}
\left[ \partial^\gamma D^\beta G(x-y,t) -\partial^\gamma D^\beta
G(x,t)\right]I_\beta u_0(y)\;dy\right|dx \label{e-AG} \\
&&\leq 
\int\!\! \int_{\bbfR^n\times \bbfR^n}|I_\beta u_0(y)|
\left| (\partial^\gamma D^\beta G)(z-y/\sqrt t,1) -(\partial^\gamma
D^\beta G)(z,1)\right|
\;dydz\nonumber 
\end{eqnarray}

From the first part of this proof,  the function
$\partial^\gamma D^\beta G(z,1)$ is continuous, hence the integrand
on the right hand side of \rf{e-AG} tends to 0 as $t\to\infty$ for
all $y,z\in \bbfR^n$. Denote 
$$
{\cal A}(z,y,t)\equiv (\partial^\gamma D^\beta G)(z-y/\sqrt t,1)
-(\partial^\gamma D^\beta G)(z,1).
$$
To apply the Lebesgue Dominated Convergence Theorem
to the integral on the right hand side of \rf{e-AG}, 
it is necessary to 
show that there exists $F\in L^1(\bbfR^n)$ independent
of $y\in \bbfR^n $ and $t\geq 1$, such that
\begin{equation}
|{\cal A}(z,y,t)|\leq F(z)\label{A-estim}
\end{equation}
for all $z,y\in \bbfR^n$ and $t\geq 1$.
Note  that 
$$
 {\cal A}(z, y,t) =\int_{R^n}
 |\xi|^{\beta}(i\xi)^{\gamma}\left[e^{-iy/\sqrt{t}}-1\right]
e^{-|\xi|^2} e^{iz\xi} d\xi.
$$
Moreover, the symbol 
$b(\xi,y,t) \equiv (1+|\xi|^2)^{\beta/2} (i\xi)^\gamma 
\left[e^{-iy/\sqrt{t}}-1\right]
e^{-|\xi|^2}$ is a $C^\infty$ function of 
$(\xi,y)\in \bbfR^n\times\bbfR^n$, and satisfies the differential
inequalities 
$$
|\partial^\alpha_\xi \partial ^\gamma_y b(\xi,y,t)|\leq C(\alpha,
\gamma,N) (1+|\xi|)^{-N-\alpha}
$$
for all multi-indices $\alpha$ and $\gamma$, all $N\in I\!\!N$, and 
$C(\alpha, \gamma,N)$ independent of $\xi,y\in\bbfR^n$ and $t\geq 1$.
By \cite[Ch. VI, Sec. 4, Prop. 1]{St}, the (inverse) Fourier
transform with respect to $\xi$ of $b(\xi,y,t)$
satisfies the estimate 
$$
|{\cal F}^{-1}_\xi b(\cdot,y,t)(z)|\leq C(1+|z|)^{-N}
$$
for all $N\in I\!\!N$, and a constant
$C=C(N)$ independent of $z,y\in\bbfR^n$ and $t\geq 1$.

Finally, the use of the measure $\mu_\beta$ from the first part of
this proof combined with standard properties of the Fourier transform
and the convolution lead to the representation
$
{\cal A} (\cdot, y,t) =\mu_\beta*{\cal F}^{-1}_\xi b(\cdot,y,t).
$
Hence, \rf{A-estim} holds true for the function 
$F(z)=C[\mu_\beta*(1+|\cdot|^{-N})](z)$ with any $N>n$.
This completes  the proof of 
Proposition~\ref{prop-L1-lin}.
\cbdu

\bigskip

We recall that, in \cite{MI}, Miyakawa  obtained the $L^1$-decay 
 of $e^{t\Delta}u_0$ 
provided  the $|x|^{\beta}$-momentum of the data is bounded.
Below, we
will
show that our assumptions is weaker than the one assumed by Miyakawa.

\begin{remark}
The $L^1$-decay of solutions to the linear heat equation formulated
in \rf{lin-L1} was proved by Miyakawa \cite{MI} under the
assumptions 
\begin{eqnarray}
u_0\in L^1(\bbfR^n), \;\; \int_{\bbfR^n} u_0(x)\;dx=0,\;\; 
\int_{\bbfR^n} |x|^\beta |u_0(x)|\; dx<\infty
\label{M:as}
\end{eqnarray}
for some $0<\beta<1$.
To show that our assumption $I_\beta u_0\in L^1(\bbfR^n)$ is weaker
than (\ref{M:as})
 it suffices to establish the inequality
\begin{equation}
\|I_\beta u_0\|_1\leq C \int_{\bbfR^n} |x|^\beta |u_0(x)|\;dx
\label{Ib-u0}
\end{equation}
valid for every $u_0$ satisfying \rf{M:as} with $\beta\in (0,1)$. Let
us sketch the proof of \rf{Ib-u0}, however, it
does not play any role in our considerations, below. It is well known
that 
$
(I_\beta u_0)(x) =C(\beta,n) \int_{\bbfR^n} |x-y|^{\beta-n}
u_0(y)\;dy
$
(in fact, this representation holds true for every $\beta\in (0,n)$).
Hence, using the assumption $\int_{\bbfR^n} u_0(y)\;dy=0$ and
changing the order of integration yield
$$
\|I_\beta u_0\|_1\leq C(\beta,n) \int_{\bbfR^n}\left(
\int_{\bbfR^n}\left| {1\over |x-y|^{n-\beta}} -{1\over
|x|^{n-\beta}}\right|\;dx\right) |u_0(y)|\;dy.
$$
Next, note that the integral with respect to $x$ in the inequality
above is finite for every $y\in\bbfR^n$, because its integrand 
$||x-y|^{\beta-n}-|x|^{\beta-n}|$ is locally integrable and 
behaves like $|x|^{\beta-1-n}$ as $|x|\to\infty$ (here, the
assumption $\beta\in (0,1)$ is crucial). Hence, by the change of
variables, it follows that
$$
\int_{\bbfR^n}\left| {1\over |x-y|^{n-\beta}} -{1\over
|x|^{n-\beta}}\right|\;dx=|y|^\beta
\int_{\bbfR^n}\left| {1\over |\omega-y/|y||^{n-\beta}} -{1\over
|\omega|^{n-\beta}}\right|\;d\omega.
$$
Since $\sup_{y\in \bbfR^n\setminus \{0\}}
\int_{\bbfR^n} ||\omega-y/|y||^{\beta-n}-|\omega|^{\beta-n}|
\;d\omega<\infty$
(the proof of this elementary fact is omitted), we obtain \rf{Ib-u0}.
\end{remark}

\bigskip

The self-similar
asymptotics of $e^{t\Delta}u_0$ 
in $L^p(\bbfR^n)$ with $p\in [2,\infty]$
can be derived under weaker assumptions on $u_0$.
This is stated in the following proposition.

\begin{proposition}\label{lem:lin:ss}
Let $\ell=\ell(\xi)$ denote a function homogeneous of degree
$\beta>0$.
 Assume that $u_0$ satisfies
\begin{equation}
\sup_{\xi\in\bbfR^n\setminus \{0\}} {\widehat u_0(\xi)\over \ell
(\xi)}<\infty
\;\;\;\;\;\;\mbox{and} \;\;\;\;\;\;
\lim_{|\xi|\to 0} {\widehat u_0(\xi)\over \ell (\xi)}=A
\label{as:u0}
\end{equation}
for some $A\in \bbfR$.
Denote by $\cal L$ the Fourier multiplier operator defined via the
formula
$
\widehat{{\cal L}v}(\xi)= \ell(\xi)\widehat v(\xi).
$
Under these assumptions, 
for every $p\in [2,\infty]$ and for every multi-index $\gamma$,
it follows
$$
t^{n(1-1/p)/2+\beta/2+ |\gamma|/2} \| \partial^\gamma e^{t\Delta} u_0 
-A  \partial^\gamma {\cal L} G(t)\|_p\to 0
\quad \mbox{as}\quad  t\to\infty.
$$
\end{proposition}

\proof
The main tool here is the Hausdorff--Young inequality
\begin{equation}
\|\widehat v\|_p\leq C\|v\|_q,\label{HY}
\end{equation}
valid for every $1\leq q\leq 2\leq p\leq \infty$ such that
$1/p+1/q=1$.
Hence,  \rf{HY}, the change of variables $\xi t^{1/2}=\omega$, and
the homogeneity of $\ell$ yield
\begin{eqnarray*}
&&\|\partial^\gamma e^{t\Delta} u_0 -A\partial^\gamma {\cal L}
G(t)\|_p^q\\
&&\hspace{1cm}\leq C\int_{\bbfR^n}
\left|(i\xi)^{\gamma}e^{-t|\xi|^2} \ell(\xi)
\left({\widehat
u_0(\xi)-A\ell(\xi)\over\ell(\xi)}\right)\right|^q\,d\xi\\
&&\hspace{1cm}= Ct^{-n/2-(\beta/2+|\gamma|/2)q}\int_{\bbfR^n}
\left|(i\omega)^{\gamma}e^{-|\omega|^2}
\ell(\omega)
\left({\widehat
u_0(\omega/t^{1/2})\over\ell(\omega/t^{1/2})}-A\right)
\right|^q\,d\xi.
\end{eqnarray*}
Now, the assumptions on $u_0$ in \rf{as:u0} allow us to  
apply the Lebesgue
Dominated Convergence Theorem in order to prove that the integral on
the right hand side tends to 0 as $t\to\infty$. 
\cbdu

\medskip

\begin{remark}
The conditions formulated in \rf{as:u0} appear in a natural way 
if  Hardy spaces are considered. 
Let us recall that a tempered distribution $v$ belongs to the Hardy 
space ${\cal H}^p$ on $\rn$ for some $0< p <\infty$ whenever
$ 
v^+
= \sup_{t>0} |(\phi_t * v)| \in L^p({\BR^n}),
$
where $\phi_t(x) = t^{-n} \phi(x/t)$ with $\phi \in {\cal S}(\rn)$
such
 that $\int_{\BR^n} \phi(x) \dd x =1$. 
We refer the reader to \cite{St} where several properties of 
Hardy spaces are derived.  
We recall that  ${\cal H}^1$  
is a Banach space strictly contained in $L^1(\bbfR^n)$ and that
$L^p(\rn) 
= {\cal H}^p$ for $p >1$ with  equivalent norms.
Suppose now  that $p\leq 1$ and $u_0\in \H^p$. It is known 
(cf. \cite[Chapter III, \S 5.4]{St}) that the Fourier transform
$\widehat u_0$ is continuous on $\bbfR^n$ and
$
|\widehat u_0(\xi)|\leq C |\xi|^{n(1/p-1)}\|u_0\|_{\H^p}
$
for all $\xi\in\bbfR^n$.
Moreover, near the origin, this can be refined to
$
\lim_{\xi\to 0} {\widehat u_0(\xi) |\xi|^{-n(1/p-1)}}=0.
$
Hence,  assumptions \rf{as:u0} are satisfied with
$\ell(\xi)=|\xi|^{\beta}$, $\beta\in (0,1)$, and $A =0$, if e.g.
$u_0\in
\H^{n/(n+\beta)}$.
\end{remark}

\bigskip


 Theorem \ref{tw:L1dec} is the main decay theorem proved in  this
 section.
 It ensures that  the $L^1$-norm of  solutions to the 
 convection-diffusion 
equation decay at the rate  $t^{-\beta/2}$ provided  their initial
data are such
that
 the corresponding solutions to the
 heat equation decay at  the same rate.

\bigskip

\noindent{\bf Proof of Theorem \ref{tw:L1dec}.}
The proof of this theorem relays on a systematic combination of  the 
integral equation \rf{duhamel} 
with  inequality \rf{CL}. 
Note that since $u_0\in L^1(\bbfR^n)\cap L^q(\bbfR^n)$,
 by \rf{CL}, it follows that
\begin{equation}
\|u(\cdot,t)\|_q^q \leq C(\|u_0\|_1,\|u_0\|_q)(1+t)^{-(n/2)(q-1)}
\label{L1-dec0}
\end{equation}
for all $t\geq 0$. Hence, computing the $L^1$-norm of \rf{duhamel},
 using the assumption on $u_0$, and \rf{L1-dec0} yield
\begin{eqnarray}
\|u(\cdot,t)\|_1&\leq& \|e^{t\Delta}u_0\|_1 +
\int_0^t\|a\cdot\nabla G(\cdot,t-\tau)\|_1\|u(\cdot,t)\|_q^q\,d\tau
\nonumber\\
&\leq& Ct^{-\beta/2} + C\int_0^t (t-\tau)^{-1/2} 
(1+\tau)^{-(n/2)(q-1)}\,d\tau \label{L1-dec1}\\
&\leq& Ct^{-\beta/2} +  
 C
\left\{
\begin{array}{l@{\quad\mbox{for}\quad}c}
t^{1/2-(n/2)(q-1)}, &q\in \left(1+{1\over n},1+{2\over n}\right);\\
t^{-1/2}\log (e+t), &q=1+{2\over n};\\
t^{-1/2}, &q>1+{2\over n}.
\end{array}
\right.\nonumber
\end{eqnarray}

For $q\geq 1+(\beta+1)/n$, estimate  \rf{L1:dec} follows immediately
from \rf{L1-dec1}, since $-1/2 <-\beta/2$ and  $1/2-(n/2)(q-1)\leq
-\beta/2$ for this
range of $q$. 

Next,  consider $ 1+1/n < q < 1+(\beta+1)/n$.
A simple calculation shows that
$\alpha =-(1/2-(n/2)(q-1))$
satisfies $0<\alpha<\beta/2$. Moreover, it follows from \rf{L1-dec1}
that 
\begin{equation}
\|u(\cdot,t)\|_1\leq C(1+t)^{-\alpha}.\label{L1-dec1a}
\end{equation}
Combining inequality  \rf{CL} with \rf{L1-dec1a}  
yields the improved decay of
the $L^q$-norm  
\begin{eqnarray}
\|u(\cdot,t)\|_q&\leq& C(1+t/2)^{-(n/2)(1-1/q)}\|u(\cdot,t/2)\|_1
\label{L1-dec2}\\
&\leq& C(1+t)^{-(n/2)(1-1/q)-\alpha}\nonumber.
\end{eqnarray}
Hence, repeating the calculations from \rf{L1-dec1}, using
\rf{L1-dec2}
instead of \rf{L1-dec0}, gives
\begin{equation}
\|u(\cdot,t)\|_1\leq Ct^{-\beta/2} +
\int_0^t(t-\tau)^{-1/2} (1+\tau)^{-(n/2)(q-1)-q\alpha}d\tau.
\label{L1-dec3}
\end{equation}
If $-(n/2)(q-1)-q\alpha\leq -1$, the integral on the right hand
side
of \rf{L1-dec3} tends to 0 as $t\to\infty$ faster than
$t^{-\beta/2}$
and this ends the proof. On the other hand, if 
$-(n/2)(q-1)-q\alpha> -1$, by the definition of $\alpha$, it
follows
from \rf{L1-dec3} that
$$
\|u(\cdot,t)\|_1\leq Ct^{-\beta/2} +Ct^{-\alpha(q+1)}.
$$
Hence, if $-\alpha(q+1)\leq -\beta/2$, the proof is complete. If, on
the contrary
 $-\alpha(q+1) > -\beta/2$,      we have the new estimate 
$$
\|u(\cdot,t)\|_1\leq C(1+t)^{-\alpha(q+1)},
$$
which we
 use as in \rf{L1-dec2} (with $\alpha$
 replaced by $\alpha(q+1)$)
to get an improved decay  of the $L^q$-norm:
$ \|u(\cdot,t)\|_q\leq 
C(1+t)^{-(n/2)(1-1/q)-\alpha(q+1)}$.
Consequently, 
 a finite number of repetitions of
the above steps yields \rf{L1:dec}.

Finally, let us prove \rf{L1:dec} for $1+1/(n+\beta)\leq q\leq 1+1/n$
under
the assumption that $\sup_{t>0}\|e^{t\Delta}u_0\|_1$ is sufficiently
small. For simplicity of notation, we put
$$
q=q^*=1+{1\over n+\beta},
$$
and we use systematically the following inequality 
(obtained from the H\"older inequality and from \rf{CL})
\begin{equation}
\|u(\cdot,t)\|_q^q\leq
\|u(\cdot,t)\|_{q^*}^{q^*}\|u(\cdot,t)\|_\infty^{q-q^*}
\leq C(\|u_0\|_\infty)\|u(\cdot,t)\|_{q^*}^{q^*}
\label{q-q*}
\end{equation}
for all $t>0$.
To proceed, we also  define the auxiliary nonnegative continuous
function 
$$
g(t)\equiv \sup_{0\leq \tau\leq t}
\left(\tau^{\beta/2}\|u(\cdot,\tau)\|_1\right)  
+\sup_{0\leq \tau\leq t}
\left(\tau^{(1/2+\beta/2)/q^*} \|u(\cdot,\tau)\|_{q^*}\right).
$$

Now, computing the $L^1$-norm of the integral equation \rf{duhamel}
and
using \rf{q-q*} yield
\begin{eqnarray}
t^{\beta/2}\|u(\cdot,t)\|_1&\leq& t^{\beta/2}\|e^{t\Delta}u_0\|_1
+Ct^{\beta/2}\int_0^t(t-\tau)^{-1/2}
\|u(\cdot,\tau)\|_{q^*}^{q^*}\;d\tau
\nonumber\\
&\leq& 
t^{\beta/2}\|e^{t\Delta}u_0\|_1\label{1-est} \\ 
&&+ \;\; g^{q^*}(t)\; C
t^{\beta/2}\int_0^t(t-\tau)^{-1/2}
\tau^{-1/2-\beta/2}\;d\tau\nonumber
\end{eqnarray}
for all $t>0$. An elementary calculation shows that the quantity
$$
t^{\beta/2}\int_0^t(t-\tau)^{-1/2}
\tau^{-1/2-\beta/2}\;d\tau
$$ 
is finite for every $t>0$ (since
$0<\beta<1$) and independent of $t$.
A similar reasoning gives
\begin{eqnarray}
\|u(\cdot,t)\|_{q^*}&\leq&
(t/2)^{-(n/2)(1-1/q^*)}\|e^{(t/2)\Delta}u_0\|_1
\label{q*-est}\\
&&+\; g^{q^*}(t)\,C\int_0^t(t-\tau)^{-(n/2)(1-1/q^*)-1/2} 
\tau^{-1/2-\beta/2}\;d\tau. \nonumber
\end{eqnarray}
Note  that $-(n/2)(1-1/q^*)-\beta/2=-(1/2+\beta/2)/q^*$.
Moreover, the quantity 
$$
t^{(1/2+\beta/2)/q^*}\int_0^t(t-\tau)^{-(n/2)(1-1/q^*)-1/2} 
\tau^{-1/2-\beta/2}\;d\tau
$$ 
is finite (since $-(n/2)(1-1/q^*)-1/2>-1$) and independent of $t$
(by
the change of variables).

Combining inequalities \rf{1-est} and \rf{q*-est} yields
\begin{equation}
g(t)\leq C_1 \sup_{0\leq \tau}\tau^{\beta/2} \|e^{t\Delta} u_0\|_1
+C_2 g^{q^*}(t)
\label{gt-est}
\end{equation}
for all $t\geq 0$ and constants $C_1$ and $C_2$ independent of $t$.

Finally, let
$$
F(y)=A+C_2y^{q^*} -y\quad \mbox{where}\quad  
A=C_1 \sup_{0\leq \tau}\tau^{\beta/2} \|e^{t\Delta} u_0\|_1
$$ 
and where
$q^*>1$. If $A>0$ is sufficiently small, there exists $y_0>0$ such
that $F(y_0)=0$ and $F(y)>0$ if $y\in [0, y_0)$. Moreover, it follows
from \rf{gt-est} that $F(g(t))\geq 0$. 
Since $g(t)$ is a nonnegative, continuous function such that
$g(0)=0$,
we deduce that $g(t)\in [0, y_0)$ for all $t\geq 0$.
This
completes the proof of Theorem \ref{tw:L1dec}.
\cbdu

\bigskip


As a consequence of  Theorem   \ref{tw:L1dec} we get  Corollary
\ref{cor-Lp0}. Actually, here we prove its
 slightly stronger version.

\begin{corollary}
\label{cor-Lp}
Under the assumptions of Theorem \ref{tw:L1dec}, for every $p\in
[1,\infty]$ and $\beta \in (0,1)$
there exists $C=C(u_0,p)$ independent of $t$   
such that
\begin{equation}
\|u(\cdot,t)\|_p\leq C(1+t)^{-(n/2)(1-1/p)-\beta/2} \label{u-Lp}
\end{equation}
for all $t>0$, and
\begin{eqnarray}
&&\hspace{-1cm}\|u(\cdot,t)-
e^{t\Delta}u_0(\cdot)\|_p\nonumber\\
& & \label{ue-Lp} \\
&\leq& C 
\left\{
\begin{array}{l@{\quad\mbox{for}\quad}c}
t^{-(n/2)(q-1/p)-(\beta q-1)/2} &q\in
\left(1+{1\over n+\beta},
{n+2\over n+\beta}\right),\\
t^{-(n/2)(1-1/p)-1/2}\log(e+ t) &q= {n+2\over n+\beta},\\
t^{-(n/2)(1-1/p)-1/2} & q>{n+2\over n+\beta}
\end{array}
\right.\nonumber
\end{eqnarray}
for all $t\geq 1$.
\end{corollary}

\proof 
Inequality \rf{u-Lp} is obtained combining \rf{CL} with \rf{L1:dec}
as
in \rf{L1-dec2} where $q$ is replaced by $p$  and $\alpha$ by
$\beta/2$.

In view of the integral equation \rf{duhamel}, to prove \rf{ue-Lp}
it
suffices to estimate the $L^p$-norm of the second term on the right
hand
side of \rf{duhamel}.  Here, split the integration range with
respect to $\tau$ into $[0,t/2]\cup [t/2,t]$ and study each term
separately as follows. Using the Young inequality for the convolution
and \rf{u-Lp} yields
\begin{eqnarray}
&&\hspace{-1cm}
\int_0^{t/2}\|a\cdot\nabla
e^{(t-\tau)\Delta}(u|u|^{q-1})(\tau)\|_p\;d\tau
\nonumber \\
&\leq& \int_0^{t/2}\|a\cdot\nabla G(\cdot,t-\tau)\|_p
\|u(\cdot,\tau)\|_q^q\;d\tau\label{ue-Lp-1}\\
&\leq& C\int_0^{t/2}(t-\tau)^{-(n/2)(1-1/p)-1/2}
(1+\tau)^{-(n/2)(q-1)-\beta q/2}\;d\tau\nonumber\\
&\leq& C
\left\{
\begin{array}{l@{\quad\mbox{for}\quad}c}
t^{-(n/2)(q-1/p)-(\beta q-1)/2} &q\in\left(1+{1\over n+\beta},
{n+2\over n+\beta}\right),\\
t^{-(n/2)(1-1/p)-1/2}\log(e+ t) &q= {n+2\over n+\beta},\\
t^{-(n/2)(1-1/p)-1/2} &q> {n+2\over n+\beta}
\end{array}
\right.\nonumber
\end{eqnarray}
for all $t>0$.

A similar calculation gives
\begin{eqnarray}
&&\hspace{-2cm}
\int_{t/2}^t\|a\cdot\nabla
e^{(t-\tau)\Delta}(u|u|^{q-1})(\tau)\|_p\;d\tau
\nonumber\\
&\leq& 
\int_{t/2}^t (t-\tau)^{-1/2}
 \|u(\cdot,\tau)\|_{pq}^q\;d\tau\label{ue-Lp-2}\\
&\leq & C\int_{t/2}^t(t-\tau)^{-1/2}
(1+\tau)^{-(n/2)(q-1/p)-\beta q/2}\;d\tau\nonumber\\
&\leq& C t^{-(n/2)(q-1/p)-(\beta q-1)/2}
\nonumber
\end{eqnarray}
for all $t>0$.

Finally, to obtain \rf{ue-Lp}, combine \rf{ue-Lp-1} and \rf{ue-Lp-2}
(note that 
$$
-(n/2(q-1/p)-(\beta q-1)/2 \leq -(n/2)(1-1/p)-1/2
$$
for 
$q\geq (n+2)/(n+\beta)$). 
\cbdu

\bigskip

{\bf Proof of Corollary \ref{cor-lin-self}.}
As pointed out in Section 2, it follows from  Corollary~\ref{cor-Lp}
that
\begin{equation}
t^{(n/2)(1-1/p)+\beta/2} \|u(\cdot,t) -e^{t\Delta}u_0(\cdot)\|_p\to 0 
\quad \mbox{as} \quad t\to\infty,
\label{u-e-0}
\end{equation}
 for each $q>1+1/(n+\beta)$   and  every $p\in [1,\infty]$. 
Thus
 \rf{u-e-0}  combined with  Propositions
\ref{prop-L1-lin}
and \ref{lem:lin:ss} yields Corollary \ref{cor-lin-self}. 
\cbdu

\medskip

A few remarks are  in order.

\begin{remark} 
If the nonlinear term in \rf{eq} has the form $\nabla\cdot f(u)$ and
the function $f$ is sufficiently regular at zero, it is possible to
improve the conclusion  of Corollary \ref{cor-lin-self} to  
$$
t^{(n/2)(1-1/p)+(\beta+|\gamma|)/2}
\|\partial^\gamma u(\cdot,t)-A\partial^\gamma 
D^\beta G(\cdot,t)\|_p\to 0
\quad \mbox{as} \quad
t\to\infty
$$
for the multi-index $\gamma$ depending on the regularity of $f$.
\end{remark}

\begin{remark}
Consider Theorem \ref{tw:L1dec} and  Corollary
\ref{cor-lin-self} 
 in the context of the
viscous Burgers equation
\begin{equation}
u_t-u_{xx}+(u^2/2)_x=0,\quad u(x,0)=u_0(x), \quad x\in\bbfR.
\label{burgers}
\end{equation}
This is  problem \rf{eq}-\rf{ini} with $n=1$, $q=2$, and $a=1/2$. 
It is well-known (cf. e.g. \cite{H,CTPL,EZ,BKW3,DuZ,DuZ2}) that the
large time behavior of solutions to this equation supplemented with
the integrable
initial condition is described by so-called {\it nonlinear diffusion
waves} (cf. \rf{auto}, above). 
If, however, it is assumed   that $u_0$ satisfies
the conditions from Proposition \ref{prop-L1-lin}
 with some $0<\beta<1$, and if, moreover,
 $\sup_{t>0}\|e^{t\Delta}u_0\|_1$ is sufficiently small,  the
 asymptotics for large $t$ of solutions to the Burgers
equation is given by the self-similar solutions $AD^\beta G(x,t)$ to
{\it the heat equation}.

For  completeness of the exposition,  we analyze problem
\rf{burgers} in more detail. 
Using the Hopf-Cole
transformation one obtains the solution of \rf{burgers} of the
following form
\begin{equation}
u(x,t) = -
{(e^{t\Delta} w_0)_x(x)\over (e^{t\Delta}w_0)(x)}
\label{cole-hopf}
\end{equation}
where as usual 
$ w_0(x) = \exp\left(- \int_{-\infty}^x u_0(y)\; dy\right)$.
Supposing
that
$u_0\in L^1(\bbfR)$, $\int_{\bbfR} u_0(x)\;dx=0$, and $u_0$ satisfies
the Miyakawa moment condition
\begin{equation}
\int_{R^n} |x|^{\beta}|u_0(x)|dx < \infty, \label{mommi}
\end{equation}
it is 
easy to show directly from the explicit formula
(\ref{cole-hopf}) that the $L^1$-norm of $u(\cdot,t)$  decays at the
rate  $t^{-\beta/2}$. Indeed, first  note that 
the denominator is the solution to the heat equation with the
datum
$w_0$ and  is bounded from below by
$\exp(-\|u_0\|_1)$. 
Thus, it is only necessary  to bound the numerator 
\begin{eqnarray*}
\left(e^{t\Delta} w_0\right)_x(x) &=&
\int_{\bbfR}G(x-y,t)(w_0)_y(y)\;dy\\
&=& - \int_{\bbfR}G(x-y,t) u_0(y)
\exp\left(-\int_{-\infty}^yu_0(z)\;dz\right)\;dy.
\end{eqnarray*}
Obviously, $(w_0)_y\in L^1(\bbfR)$ and $\int_{\bbfR} |y|^\beta
|(w_0)_y(y)|\;dy<\infty$, since $u_0$ has these properties and
$\exp\left(-\int_{-\infty}^yu_0(z)\;dz\right)$ is a bounded function.
Let us skip an easy proof that $\int_{\bbfR} (w_0)_y(y)\;dy=0$. 
Consequently, $(w_0)_y$ satisfies the Miyakawa conditions, so the
$L^1$-norms of solutions to \rf{burgers} decay  with the rate
$t^{-\beta/2}$. Finally, repeating the calculations from the proof of
Corollary \ref{cor-Lp} yields that the large time behavior of
solutions to
\rf{burgers} is described by $AD^\beta G(\cdot,t)$. Note that here 
non-smallness  assumptions on $u_0$  have been imposed unlike
 it was done in
Theorem \ref{tw:L1dec} in the case $1+1/(n+\beta)\leq q \leq 1+1/n$.
This example suggests that such an assumption in Theorem
\ref{tw:L1dec} is not necessary, however, the proof of a stronger
version requires new ideas.
\end{remark}

\begin{remark}
In this paper,  we limit ourselves to the case $\beta\in
(0,1)$. We expect a completely different large time behavior of
solutions to \rf{eq}-\rf{ini} in $\beta\geq 1$ for the following
reason. Suppose that 
\begin{equation}
u_0\in L^1(\bbfR^n, (1+|x|)\;dx) \quad \mbox{and} \quad 
\int_{\bbfR^n} u_0(x)\; dx=0.
\label{r:u0}
\end{equation}
It is proved in \cite{DZ92} 
that 
$
\|e^{t\Delta}u_0\|_1\leq Ct^{-1/2}\|u_0\|_{L^1(\bbfR^n, |x|\;dx)}
$
for all $t>0$ and a constant $C$; moreover,   
$$
t^{1/2}\left\|e^{t\Delta}u_0-\int_{\bbfR^n} xu_0(x)\;dx \cdot 
\nabla G(x,t)\right\|_1\to 0 \quad \mbox{as}\quad t\to\infty.
$$
Now, using   the second order asymptotic expansion by Zuazua
\cite{Z} (cf. also \cite{BKW2} for analogous results with more
general diffusion operators and less regular initial conditions) of
solutions to \rf{eq}-\rf{ini} with $q>1+2/n$, we obtain that the
quantity 
$$
t^{1/2}
\left\|u(\cdot,t)
-\left(\int_{\bbfR^n} xu_0(x)\;dx 
-a\int_0^\infty\int_{\bbfR^n}(u|u|^{q-1})(x,\tau)\;dxd\tau\right)
\cdot \nabla G(x,t)\right\|_1
$$
tends to 0 
as $t\to\infty$. This asymptotic result shows that the large time
behavior of solutions with the initial data satisfying \rf{r:u0}
can be classified as weakly nonlinear in the sense of Zuazua
\cite{Z}. Here, however, the first term of the asymptotics comes
linearly from the
heat kernel, but has a nonlinear dependence on the solution through a
multiplicative factor (as noted by Zuazua in \cite{Z}, it is an open
question if this factor is different from zero).
Hence,
assuming that $\|e^{t\Delta}u_0\|_1\leq Ct^{-\beta/2}$ for some
$\beta\geq 1$ one should expect asymptotic expansions of solutions
completely different from that in Corollary \ref{cor-lin-self},
specifically of the form just described.
\end{remark}


\setcounter{equation}{0}
\section{Nonlinear asymptotics}

The following two lemmata give the crucial  steps to yield the
necessary estimates of 
the integral equation \rf{duhamel}.

\begin{lemma}\label{lem:nabla:e}
Let  $a\in\bbfR^n$ be a fixed constant vector. There exists 
a constant $C>0$ such that for every $w\in L^1(\bbfR^n)$ we have
\begin{equation}
\|a\cdot\nabla e^{t\Delta}w\|_{\B1b}\leq Ct^{(\beta-1)/2}\|w\|_1
\end{equation}
for all $t>0$.
\end{lemma}

\proof
Using the definition of the norm in
$\B1b$ and properties of the heat semigroup yields
\begin{eqnarray*}
\|a\cdot\nabla e^{t\Delta}w\|_{\B1b}&=&
\sup_{s>0}s^{\beta/2} \|e^{s\Delta}a\cdot\nabla e^{t\Delta} w\|_1\\
&=& \sup_{s>0} s^{\beta/2} \|a\cdot\nabla e^{(t+s)\Delta}w\|_1\\
&\leq& C\|w\|_1 \sup_{s>0} s^{\beta/2}(t+s)^{-1/2}
\end{eqnarray*}
for all $t>0$. Now, a direct calculation shows that
$\sup_{s>0}s^{\beta/2}(t+s)^{-1/2}=C(\beta) t^{(\beta-1)/2}$ with
$C(\beta)$
independent of $t$.
\cbdu

\begin{lemma}\label{lem:eB}
Assume that $v\in \B1b$. Then for each $p\in [1,\infty]$ there exists
a constant $C>0$ such that 
$$
\|e^{t\Delta}v\|_p\leq Ct^{-(n/2)(1-1/p) -\beta/2}\|v\|_\B1b
$$
for all $t>0$.
\end{lemma}

\proof
Standard properties of the heat semigroup $e^{t\Delta}$
and the definition of the norm in $\B1b$ give
$$
\|e^{t\Delta}v\|_p\leq
C(t/2)^{-(n/2)(1-1/p)}\|e^{(t/2)\Delta}v\|_{1}\leq 
Ct^{-(n/2)(1-1/p) -\beta/2}\|v\|_\B1b.
$$
for all $t>0$ and a constant $C$.
\cbdu

\bigskip


We are ready to prove the existence Theorem 
\ref{exist:global} in the critical case $q^*$.

\medskip

\noindent {\bf Proof of Theorem \ref{exist:global}.}
Our reasoning is similar to that in \cite{C,C2,CP,K}. Moreover, the
calculations below resemble those in the proof of Theorem
\ref{tw:L1dec} with $1+1/(n+\beta)\leq q\leq 1+1/n$,
thus we shall be brief in details. 
Recall that in this section we consider
$$
q=q^*=1+{1\over n+\beta}
$$
which is equivalent to 
$$
{n\over 2}\left(1-{1\over q}\right) +{\beta\over 2} =
{1\over q} \left({1\over 2} +{\beta\over 2}\right).
$$

Equip the space $\X$ with the norm 
$$
\|u\|_{\X}=\max \{\sup_{t>0}\|u(t)\|_{\B1b},\;\; 
\sup_{t>0}t^{(n/2)(1-1/q)+\beta/2}\|u(t)\|_q\}.
$$
We will show that the nonlinear operator 
\begin{equation}
\Nop(u)(t)\equiv e^{t\Delta}u_0 -\int_0^t a\cdot \nabla
e^{(t-\tau)\Delta} (u|u|^{q-1})(\tau)\,d\tau \label{Nop} 
\end{equation}
is a  contraction on the box
$$
B_{R,\varepsilon}=\{u\in \X \;:\; \|u(t)\|_{\B1b}\leq
R\;\;\mbox{and}\;\;  
\sup_{t>0}t^{(n/2)(1-1/q)+\beta/2}\|u(t)\|_q\leq 2\varepsilon\}
$$
for sufficiently large $R>0$ and a suitably small $\varepsilon>0$.
This will be guaranteed  provided  the following
estimates can be shown to hold
\begin{eqnarray}
\|\Nop(u)(t)\|_\B1b&\leq&
\|u_0\|_\B1b+C\varepsilon^q,\label{es:1}\\
t^{(n/2)(1-1/q)+\beta/2}\|\Nop(u)(t)\|_q&\leq& 
C\|u_0\|_\B1b+C\varepsilon^q,\label{es:2}
\end{eqnarray}
and
\begin{eqnarray}
&&\|\Nop(u)(t)-\Nop(v)(t)\|_\B1b\label{es:3}\\
&&\hspace{1cm}\leq
C\varepsilon^{q-1}\sup_{t>0}t^{(n/2)(1-1/q)+\beta/2}\|u(\cdot,
t)-v(\cdot,t)\|_q\nonumber\\
&&t^{(n/2)(1-1/q)+\beta/2}\|\Nop(u)(t)-\Nop(v)(t)\|_q\label{es:4}\\
&&\hspace{1cm}\leq
C\varepsilon^{q-1}\sup_{t>0}t^{(n/2)(1-1/q)+\beta/2}\|u(\cdot,
t)-v(\cdot,t)\|_q.\nonumber
\end{eqnarray}
with constants $C$ independent of $u$ and $t$.

For the proof of \rf{es:1} observe that 
$\|e^{t\Delta} u_0\|_\B1b\leq \| u_0\|_\B1b$. Hence computing the
$\B1b$-norm of \rf{Nop} for $u\in B_{R,\varepsilon}$ and applying
Lemma \ref{lem:nabla:e} we obtain
\begin{eqnarray*}
\|\Nop(u)(t)\|_\B1b&\leq &\|e^{t\Delta}u_0\|_\B1b +\int_0^t \|a\cdot
\nabla e^{(t-\tau)\Delta} (u|u|^{q-1})(\tau)\|_\B1b\,d\tau\\
&\leq &\|u_0\|_\B1b +C\int_0^t (t-\tau)^{(\beta-1)/2}
\|u(\tau)\|_q^q\,d\tau\\
&\leq &\|u_0\|_\B1b +C\varepsilon^q \int_0^t (t-\tau)^{(\beta-1)/2} 
\tau^{-(n/2)(q-1)-\beta q/2}\,d\tau. 
\end{eqnarray*}
Note now that the assumptions $\beta\in (0,1)$ and $q=1+1/(n+\beta)$
guarantee that the integral on the right hand side is finite for any
$t>0$. Moreover, since $(\beta-1)/2-n(q-1)/2-\beta q/2+1=0$, it
follows that this integral is independent of $t$. Hence, estimate
\rf{es:1} holds true.

The proof of \rf{es:2} is similar. It involves Lemma \ref{lem:eB} as
follows
\begin{eqnarray}
\|\Nop(u)(t)\|_q&\leq &\|e^{t\Delta}u_0\|_q +\int_0^t \|a\cdot \nabla
e^{(t-\tau)\Delta} (u|u|^{q-1})(\tau)\|_q\,d\tau\nonumber\\
&\leq& Ct^{-(n/2)(1-1/q)-\beta/2}\|u_0\|_\B1b\label{Nop:q}\\
&&\mbox{} +C\varepsilon^q \int_0^t
(t-\tau)^{-(n/2)(1-1/q)-1/2}\tau^{-(n/2)(q-1)-\beta q/2}\,d\tau.
\nonumber
\end{eqnarray}
In this case, the conditions on $\beta, q$ imply  again that the
integral on the right hand side is finite for every $t>0$. In fact,
by a change of variables, it equals $Ct^{-(n/2)(1-1/p)-\beta/2}$ for
a
constant $C>0$. Hence \rf{es:2} is proved.

The proofs of \rf{es:3} and \rf{es:4} are completely analogous. The
only difference consists in using elementary inequality
\begin{equation}
\left\|u|u|^{q-1}-v|v|^{q-1}\right\|_1\leq C
\|u-v\|_q\left(\|u\|_q^{q-1}+\|v\|_q^{q-1}\right)\label{elem}
\end{equation}
valid for all $u,v\in L^q(\bbfR^n)$.

Finally, it follows from \rf{es:1}--\rf{es:4} that 
$\Nop : B_{R,\varepsilon}\to B_{R,\varepsilon}$ is a 
contraction for $R>2\|u_0\|_\B1b$ and a suitably small
$\varepsilon>0$. 
Hence the sequence defined as $u_0(t)=e^{t\Delta}u_0$ and
$u_{n+1}(t)=\Nop(u_n(t))$ converges to a unique (in
$B_{R,\varepsilon}$) global-in-time solution to \rf{eq}-\rf{ini}
provided $u_0(t)\in B_{R,\varepsilon}$, i.e. $\|u_0\|_\B1b
$ is sufficiently small (cf. Lemma \ref{lem:eB}). 
\cbdu

\bigskip

The proof of Theorem \ref{asymp} requires the following result from
\cite[Lemma~6.1]{K}.

\begin{lemma}\label{f:g}
Let  $w\in L^1(0,1)$, $w\geq 0$, and $\int_0^1w(x)\;dx<1$.
Assume that $f$ and $g$ are two nonnegative, bounded functions such
that 
\begin{equation}
f(t)\leq g(t)+\int_0^1 w(\tau) f(\tau t)\;d\tau.\label{L:0}
\end{equation}
Then $\lim_{t\to\infty}g(t)=0$ implies $\lim_{t\to\infty} f(t)=0$.
\cbdu
\end{lemma}


The next task is to prove   the stability Theorem \ref{asymp}.

\noindent{\bf Proof of Theorem \ref{asymp}.}
The subtraction of  equation \rf{duhamel} for $v$ 
from the
analogous expression for $u$ leads to   the
following identity
\begin{eqnarray}
u(t)-v(t)&=& e^{t\Delta}(u_0-v_0)\nonumber \\
&&-\; \int_0^t a\cdot\nabla e^{(t-\tau)\Delta}
\left(u|u|^{q-1}-v|v|^{q-1}\right)(\tau)\;d\tau.
\label{u-v}
\end{eqnarray}
Repeating the reasoning from the proof of \rf{es:2} involving
inequality \rf{elem} gives
\begin{eqnarray}
&&\hspace{-1.5cm}\|u(\cdot,t)-v(\cdot,t)\|_q\nonumber\\
&\leq&
Ct^{-(n/2)(1-1/q)-\beta/2} \left( (t/2)^{\beta/2} 
\|e^{(t/2)\Delta}(u_0-v_0)\|_1\right)\label{as:2}\\
& &+ 
C\int_0^t(t-\tau)^{-(n/2)(1-1/q)-1/2}\|u(\cdot,\tau)-v(\cdot,\tau)\|_
q
\nonumber\\
& &\hspace{3cm}
\times\left(\|u(\cdot,\tau)\|_q^{q-1}+\|v(\cdot,\tau)\|_q^{q-1}
\right)\,d\tau.\nonumber
\end{eqnarray}
By Theorem \ref{exist:global}, the both quantities
$$ 
\sup_{t>0} t^{(n/2)(1-1/q)+\beta/2}\|u(\cdot,t)\|_q
\;\;\;\mbox{and}\;\;\; 
 \sup_{t>0} t^{(n/2)(1-1/q)+\beta/2}\|v(\cdot,t)\|_q
$$
 are bounded by
$2\varepsilon$. 
Hence, multiplying \rf{as:2}  by $t^{(n/2)(1-1/q)+\beta/2}$,
putting 
\begin{equation}
f(t)=t^{(n/2)(1-1/q)+\beta/2}\|u(\cdot,t)-v(\cdot, t)\|_q,
\label{def:f}
\end{equation}
and changing variable $\tau=ts$, we get 
\begin{eqnarray}
f(t)&\leq& C(t/2)^{\beta/2}\|e^{(t/2)\Delta}(u_0-v_0)\|_1
\nonumber\\
&& + 2C\varepsilon^{q-1}
\int_0^1(1-s)^{-(n/2)(1-1/q)-1/2}s^{-(n/2)(q-1)-\beta
q/2}f(ts)\;ds.\label{as:3}
\end{eqnarray}
Since $(1-s)^{-(n/2)(1-1/q)-1/2}s^{-(n/2)(q-1)-\beta q/2}\in
L^1(0,1)$ 
(cf. comments following inequalities \rf{Nop:q}),
we may apply Lemma~\ref{f:g} obtaining
$f(t)\to 0$ as $t\to\infty$ for sufficiently small $\varepsilon>0$.
This proves  \rf{e:utvt} for $p=q$.

Next, we prove  \rf{e:utvt} for $p=1$. Computing the $L^1$-norm of
\rf{u-v}
and repeating the calculations from \rf{as:2} and \rf{as:3} yield
\begin{eqnarray*}
t^{\beta/2}\|u(\cdot,t)-v(\cdot,t)\|_1
&\leq& t^{\beta/2}\|e^{t\Delta}(u_0-v_0)\|_1 \\
&& + \; C
\int_0^1(1-s)^{-1/2}s^{-(n/2)(q-1)-\beta q/2}f(ts)\;ds,
\end{eqnarray*}
where $f$, defined in \rf{def:f}, is a bounded function satisfying
$\lim_{t\to\infty} f(t)=0$, by the 
first part of this proof. Hence \rf{e:u0v0} and the Lebesgue
Dominated Convergence
Theorem  give
\begin{equation}
\lim_{t\to\infty} t^{\beta/2} \|u(\cdot,t)-v(\cdot,t)\|_1 =0.
\label{L1:as}
\end{equation}

The next stage of the proof deals with \rf{e:utvt} for all $p\in 
(1,\infty)$. The calculations from \rf{L1-dec2} show that
$\|u(\cdot,t)\|_\infty$ and $\|v(\cdot,t)\|_\infty$ can be both
bounded by 
$Ct^{-n/2-\beta/2}$ for  all $t>0$ and a 
constant $C$ independent of $t$. Hence, by the H\"older inequality
and 
\rf{L1:as} it follows that
\begin{eqnarray*} 
\|u(\cdot,t)-v(\cdot,t)\|_p&\leq& 
C\|u(\cdot,t)-v(\cdot,t)\|_1^{1/p}\\
&& \times \left(\|u(\cdot,t)\|_\infty^{1-1/p} + 
\|v(\cdot,t)\|_\infty^{1-1/p} \right)\\
&=& o\left(t^{-(n/2)(1-1/p)-\beta/2}\right) \quad \mbox{as} \quad
t\to\infty,
\end{eqnarray*}
where we  used the following
inequality
\begin{equation}
\left|g|g|^{q-1}-h|h|^{q-1}\right|\leq
{q\over 2}|g-h|\left(|g|^{q-1}+|h|^{q-1}\right)\label{el-gh}
\end{equation}
valid for all $g,h\in\bbfR$ and  $q>1$.

Finally, the proof of \rf{e:utvt} for $p=\infty$ involves equation
\rf{u-v} and \rf{e:utvt} proved already for all $p\in [1,\infty)$. 
Standard $L^p-L^q$ estimates of the of the heat semigroup
imply that 
\begin{eqnarray*}
t^{n/2+\beta/2} \|e^{t\Delta}(u_0-v_0)\|_\infty &\leq& 
Ct^{n/2+\beta/2} (t/2)^{-n/2} \|e^{(t/2)\Delta}(u_0-v_0)\|_1\\
&=&C (t/2)^{\beta/2} \|e^{(t/2)\Delta}(u_0-v_0)\|_1 \to 0
\end{eqnarray*}
as $t\to\infty$ by assumption \rf{e:u0v0}.

To study the second term on the right hand side of \rf{u-v},
the integration range with respect to $\tau$  
is decomposed into $[0,t]=[0,t/2]\cup [t/2,t]$.

Combining  inequality \rf{el-gh} with  estimates of the heat
semi-group  
and the H\"older inequality yields
\begin{eqnarray}
&&\left\|a\cdot\nabla e^{(t-\tau)\Delta}
\left(u|u|^{q-1}-v|v|^{q-1}\right)(\tau)\right\|_\infty\nonumber\\
&&\hspace{1cm}\leq
C(t-\tau)^{-n/2-1/2} \|u(\tau)-v(\tau)\|_1
\left(\|u(\tau)\|_\infty^{q-1}+\|v(\tau)\|_\infty^{q-1}\right)
\label{u-v-1}\\
&&\hspace{1cm}\leq
C(t-\tau)^{-n/2-1/2} \tau^{-\beta/2
-(n+\beta)(q-1)/2}f_1(\tau),\nonumber
\end{eqnarray}
where $C$ is independent of $t$ and $\tau$, and 
$f_1(\tau)=\tau^{\beta/2}\|u(\tau)-v(\tau)\|_1$ is the bounded
function which tends to 0 as $t\to\infty$ by \rf{e:utvt} for $p=1$.

Moreover, choosing $1/r+1/z=1$, similar calculations lead to 
\begin{eqnarray}
&&\left\|a\cdot\nabla e^{(t-\tau)\Delta}
\left(u|u|^{q-1}-v|v|^{q-1}\right)(\tau)\right\|_\infty\nonumber\\
&&\hspace{1cm}\leq
C(t-\tau)^{-(n/2)(1-1/z)-1/2} \tau^{-(n/2)(1-1/r)-\beta/2
-(n+\beta)(q-1)/2}f_r(\tau)\label{u-v-rz}
\end{eqnarray}
where 
$f_r(\tau)=\tau^{(n/2)(1-1/r)+\beta/2}\|u(\tau)-v(\tau)\|_r$ also
tends to 0 as $t\to\infty$ by \rf{e:utvt}.
Hence, by the change of variables $\tau=ts$,
it follows from \rf{u-v-1} that
\begin{eqnarray*}
&&\hspace{-1cm}
\int_0^{t/2} 
\left\|a\cdot\nabla e^{(t-\tau)\Delta}
\left(u|u|^{q-1}-v|v|^{q-1}\right)(\tau)\right\|_\infty\;d\tau\\
&\leq&
Ct^{-n/2-\beta/2}
\int_0^{1/2} (1-s)^{-n/2-1/2}
s^{-\beta/2-(n+\beta)(q-1)/2}f_1(st)\;ds.
\end{eqnarray*}
The integral on the right hand side is finite (recall that
$q=1+1/(n+\beta)$), because 
$$
-{\beta\over 2} -{(n+\beta)(q-1)\over 2} =- {\beta+1\over 2} >-1
\quad \mbox{for $\beta\in (0,1)$}.
$$
This integral tends to 0  as $t\to\infty$ by the Lebesgue
Dominated Convergence Theorem.

The case of the integral $\int_{t/2}^t ...\;d\tau$
involves inequality \rf{u-v-rz} with $z>1$ 
chosen such that $-(n/2)(1-1/z)-1/2>-1$. 
The proof here is analogous as in  the last case and as such will be
omitted.
This completes the proof of Theorem~\ref{asymp}.

\cbdu

\bigskip


\setcounter{equation}{0}
\section{Balance case: self-similar solutions}

In this section, we continue  our analysis on the asymptotic
behaviour of solutions of (\ref{eq}) when $q$ is 
 the critical exponent $q=q^*= 1 +1/(n+\beta)$. Here, we would like
 to explain how Theorem \ref{exist:global} ensures the existence of a
 new class of self-similar solutions to \rf{eq} and how Theorem
 \ref{asymp}
 shows  that there is a large
 class of solutions whose asymptotic behaviour  in $L^p(\bbfR^n)$
 corresponds to  self-similar solutions.


Elementary calculations show that
 if $u(x,t)$ is a solution to the equation
 \begin{equation}
u_t-\Delta u+a\cdot \nabla (u|u|^{1/(n+\beta)})=0,\label{eq1}
\end{equation}
 then
so is $\lambda^{n+\beta} u(\lambda x,\lambda^2 t)$ 
for every $\lambda>0$.
Self-similar solutions should satisfy the equality
$u(x,t)= \lambda^{n+\beta} u(\lambda x,\lambda^2 t)$, hence choosing
$\lambda = \lambda(t)=1/\sqrt{t}$
yields a self-similar form
\begin{equation}
u(x,t)= t^{-{n+\beta\over 2}}U\left({x\over  \sqrt{t}}\right),
\label{U-self}
\end{equation}
where $U(x)= u(x,1)$, $x\in\bbfR^n$, and $t>0$. 
Substituting $u(x,t)$ defined in \rf{U-self} to equation \rf{eq1} we
shows
the function $U=U(x)$
satisfies the elliptic equation
\begin{equation}
-\Delta U - \frac{1}{2}x\cdot \nabla U = \frac{n+ \beta}{2}U 
+ a \cdot\nabla(U|U|^{1/(n+\beta)}) =0. \label{el}
\end{equation}
We believe that one can obtain solutions to \rf{el} using ideas
similar to those developed  in  \cite{AEZ}. In that paper, 
Aguirre, Escobedo and Zuazua establish a~priori
 estimates and existence of solutions to the system
\begin{equation}
-\Delta f - \frac{1}{2}x\cdot \nabla f = \frac{n}{2}f 
+ a \cdot\nabla \Psi(f) =0. \label{elaez}
\end{equation}
The main difference between our case and  (\ref{elaez}) is that
their 
coefficient for $f$ is exactly 
$n/2$ which is the first eigenvalue of
 $L= -\Delta f - \frac{1}{2}x\cdot \nabla f$.

In our paper, however,  we propose a completely different 
construction of self-similar solutions,
based on the Cannone method \cite{C2}.
Let us formulate this result.

\begin{theorem}\label{self-exist}
Assume that $u_0\in \B1b$ is a homogeneous distribution of degree
$-n-\beta$. Under the assumptions of Theorem \ref{exist:global}, the
constructed solution to \rf{eq}-\rf{ini} is self-similar; hence, of
the form \rf{U-self}. 
\end{theorem}

The proof of this theorem follows  the standard reasoning  (cf. e.g.
\cite[Section 3]{C2} and  \cite{C,CP,K}) and is based on the
uniqueness result
from Theorem \ref{exist:global}. Let us skip other details. Here, we
only
mention that 
 the fractional derivative of
order $\beta$ of the Dirac delta
$D^\beta \delta_0$ belongs to $\B1b$.
Indeed, this follows from the definitions of $e^{t\Delta}$
and $\delta_0$, since
$e^{t\Delta}D^{\beta}\delta_0 = D^\beta G(\cdot,t)$
(cf. the proof of Proposition \ref{prop-L1-lin}).
Hence, the self-similar form of $D^\beta G(x,t)$
(see \rf{DbG-self}) yields
$
\|e^{t\Delta}D^\beta \delta_0\|_1=t^{-\beta/2} \|D^\beta
G(\cdot,1)\|_1.
$
Finally, note that the tempered distribution $D^\beta \delta_0$ is
homogeneous of degree $-n-\beta$. Consequently, Theorem
\ref{self-exist} implies that every solution to \rf{eq}-\rf{ini}
corresponding to $A D^\beta \delta_0$ with sufficiently small $|A|$
is self-similar.

Now, let
$t^{-(n+\beta)/2}U_A(x/\sqrt t)$ denote the self-similar solution
corresponding to 
the initial datum $u_0=AD^\beta\delta_0$ for some $A\in\bbfR$.
In the following theorem, we show that $U_A$  describes 
the asymptotic behavior of a large class of solutions to
\rf{eq}-\rf{ini}.

\begin{theorem}\label{bs}
Let the assumptions from Theorem \ref{asymp} hold true.
Assume that $v$ is the solution of \rf{eq}-\rf{ini}
constructed in Theorem \ref{exist:global}
corresponding to the initial data $v_0\in\B1b$.
Let $t^{-(n+\beta)/2}U_A(x/\sqrt t)$ be the self-similar solution
corresponding to 
the initial datum $u_0=AD^\beta\delta_0$ for sufficiently small
$|A|$.
Suppose that
\begin{equation}
\lim_{t\to\infty} t^{\beta/2}\|e^{t\Delta} v_0-AD^\beta
G(\cdot,t)\|_{1}=0.\label{self-lin-as}
\end{equation}

Choosing   $\varepsilon>0$ in Theorem \ref{exist:global}
sufficiently small,  we have
$$
\lim_{t\to\infty}t^{(n/2)(1-1/p)+\beta/2}\|v(\cdot,t)-
t^{-(n+\beta)/2}U_A(\cdot/\sqrt t)\|_p= 0
$$
for every $p\in [1,\infty]$.
\end{theorem}

This theorem is  a direct corollary of Theorem \ref{asymp}.
Recall only that, by Proposition \ref{prop-L1-lin}, the limit
relation in
\rf{self-lin-as} holds true if, in
particular, $I_\beta v_0\in L^1(\bbfR^n)$. In this case,
$A=\int_{\bbfR^n} I_\beta v_0(x)\;dx$.

Let us compare Theorem \ref{bs} with its counterpart proved by
Escobedo and Zuazua in \cite{EZ}, and recalled already in
Introduction, formula \rf{auto}. When $\int u_0=M\neq 0$ and
$q=1+1/n$, equation \rf{eq} has a one-parameter family of
self-similar solutions parameterized by $M$. Moreover, $U_M$ describes
the large time asymptotics of all solutions with mass $M$. Note that
for every $u_0\in L^1(\bbfR^n)$, the condition $\int u_0=M$ is
equivalent to 
$$
\|e^{t\Delta} u_0-MG(\cdot, t)\|_1\to 0 \quad \mbox{as}\quad
t\to\infty.
$$
In our case, when $M=0$, the set of self-similar solutions to \rf{eq}
with $q=1+1/(n+\beta)$ is more complicated, however, relation
\rf{self-lin-as} (or, more generally, \rf{e:u0v0}) still allows us to identify
solutions to \rf{eq}-\rf{ini} with the given self-similar large time
behavior.


\setcounter{equation}{0}
\section{Conclusions}

The ideas developed in this paper can be applied to other types of equations. As 
the first example, let us look at the Navier-Stokes equations for the 
incompressible fluid
\begin{eqnarray*} u_t-\Delta u + u\cdot\nabla u +\nabla p 
&=& 0, \\ \mbox{div} \,u &= &  0 ,\\ u(\cdot,0) &=& u_0.  
\end{eqnarray*} 
It well-known (see 
e.g. \cite{MI}) that any integrable solenoidal 
smooth vector field $u_0$ (i.e. $ \mbox{div} \,u_0\equiv 0$) satisfies 
$\int u_0=0$. This fact motivated  Miyakawa to study in \cite{MI} the $L^2$-decay of solutions 
to the Navier-Stokes system endowed with integrable initial conditions    satisfying 
$\|e^{t\Delta}u_0\|_1\leq Ct^{-\beta/2}$ for some $0<\beta<1$, a constant $C$, and all $t>0$.
We believe that our methods will offer some improvements to the Miyakawa results.

We also expect that it will be possible to improve asymptotic expansions of solutions to the
Korteweg-de Vries-Burgers equation as well as to the Benjamin-Bona-Mahony-Burgers equation 
obtained recently in \cite{K3,K4}.  Some preliminary progress in this direction was already done by 
M. Mei in \cite{Mei}.


\bigskip
{\bf Acknowledgements.}
The authors would like to express their thanks to the anonymous
referee for
several very  helpful comments and suggestions. As a result,  the
authors
were able to
write an improved and more organized version of their paper.
Part of this research was done while G.K. was invited to Department
of Mathematics, University of California, Santa Cruz, in September
1999.  During the preparation of
the paper, G.K. was partially supported by the Foundation for Polish
Science.
Grant support from  
KBN 0050/P03/2000/18 is also gratefully  acknowledged.

\end{document}